\documentclass[oneside]{amsart}

\usepackage[scale=0.725]{geometry}

\usepackage[T1]{fontenc}
\usepackage[utf8]{inputenc}
\usepackage{lmodern}                    
\usepackage[english]{babel} 

\usepackage{csquotes}
\usepackage[style=numeric,doi=true,url=true,maxbibnames=5,backend=biber,giveninits,isbn=false]{biblatex}
\bibliography{testflow.bib}

\usepackage{amsmath, amssymb, amsthm}
\usepackage{mathtools, commath}
\usepackage{stmaryrd}
\usepackage{booktabs}
\usepackage[dvipsnames]{xcolor}
  \definecolor{darkred}{RGB}{139,0,0}
  \definecolor{mediumblue}{RGB}{0,0,205}
  \definecolor{forestgreen}{RGB}{34,139,34}
\usepackage[colorlinks=true]{hyperref}
\hypersetup{linkcolor=darkred, urlcolor=forestgreen, citecolor=mediumblue}
\usepackage{cleveref}

\usepackage[linesnumbered, noline, noend]{algorithm2e}
\SetAlFnt{\small}
\SetAlgoNlRelativeSize{-3}
\SetNlSty{}{}{}
\DontPrintSemicolon

\newcommand{\bo}{{\partial\Omega}}
\newcommand{\bog}{\Gamma}
\newcommand{\sdiv}{{\nabla\cdot}}
\DeclareMathOperator{\diver}{div}
\newcommand{\du}{{\mathcal D}(u)}
\newcommand{\dv}{{\mathcal D}(v)}

\newcommand{\famlyT}{{\mathcal T}_h}
\newcommand{\THk}[1]{\texttt{TH}_{#1}}
\newcommand{\THiso}[1]{\texttt{TH}_{#1}^\text{iso}}

\newtheorem{theorem}{Theorem}[section]

\theoremstyle{definition}

\newtheorem{algo}[theorem]{Algorithm}

\usepackage{todonotes}

\begin{document}

\title{Benchmark stress tests for flow past a cylinder at higher Reynolds numbers using EMAC}

\author[H. v. Wahl]{Henry von Wahl}
\address{Institut für Mathematik, Friedrich-Schiller-Universität Jena, Ernst-Abber-Platz 2, 07743 Jena, Germany}
\email{henry.von.wahl@uni-jena.de}

\author[L. G. Rebholz]{Leo G. Rebholz}
\address{School of Mathematical and Statistical Sciences, Clemson University, Clemson, South Carolina 29634, USA}
\email{rebholz@clemson.edu}

\author[L. R. Scott]{L. Ridgway Scott}
\address{Department of Mathematics, University of Chicago, Chicago, Illinois 60637, USA}
\email{ridg@cs.uchicago.edu}

\date{\today}

\subjclass[2010]{35Q30, 65M60}
\keywords{Drag coefficient, Chaotic flow, Strouhal period, EMAC, Taylor-Hood}

\begin{abstract}
We consider a test problem for Navier-Stokes solvers based on the flow around a cylinder at Reynolds numbers 500 and 1000, where the solution is observed to be periodic when the problem is sufficiently resolved. Computing the resulting flow is a challenge, even for exactly divergence-free discretization methods, when the scheme does not include sufficient numerical dissipation. We examine the performance of the energy, momentum and angular momentum conserving (EMAC) formulation of the Navier-Stokes equations. This incorporates more physical conservation into the finite element method even when the numerical solution is not exactly divergence-free. Consequently, it has a chance to outperform standard methods, especially for long-time simulations. We find that for lowest-order Taylor-Hood elements, EMAC outperforms the standard convective formulations. However, for higher-order elements, EMAC can become unstable on under-resolved meshes.
\end{abstract}

\maketitle

\section{Introduction}

Recently, a test problem for Navier-Stokes solvers has been proposed \cite{lrsBIBkl}.
It is a classic fluid problem, the unsteady flow around a cylinder. In contrast to the well-established benchmark \cite{schafer1996benchmark}, the recently proposed version considers higher Reynolds numbers, is posed in a larger domain, and appears to require longer simulations in time to establish the quantities of interest. Consequently, this may be more appropriate than similar test problems in the literature for higher Reynolds numbers and longer-time simulations.

The test problem was computed in \cite{lrsBIBkl,lrsBIBki} using a variety of low and higher-order finite element schemes. These consisted of Taylor-Hood elements with and without grad-div stabilization, Scott-Vogelius elements, and $H(\diver)$-conforming methods based on a mass-conserving mixed stress (MCS) formulation. In particular, it was established in \cite{lrsBIBkl} that the drag and lift of the resulting flow is periodic for $Re\leq 1000$. However, most methods required the discretization to resolve all scales before periodic behavior was observed. Only the MCS method, in conjunction with an upwind formulation of the convective term, led to the same periodic drag and lift results on coarser meshes. For increasing Reynolds numbers, periodicity was observed to be lost. Hence, we will focus on Reynolds numbers 500 and 1000 here, for which comparison values are available.

Here, we consider the conforming Taylor-Hood elements without grad-div stabilization together with the EMAC \cite{CHOR17,olshanskii2024local} variational correction scheme. The objective of EMAC is to preserve several invariants of the flow, and its name is taken from its energy, momentum, and angular momentum conservation properties. We report on these invariants for the cylinder test problem. Note that grad-div stabilization could also be included with EMAC, as was considered, e.g., in \cite{LLJ25} and \cite{CHOR17} for some computations. Taking these additional physics into account, the EMAC scheme has the potential to outperform standard schemes, particularly when long timescales are of interest~\cite{OR20}.

This paper is structured as follows. In \Cref{sec.equations}, we discuss the equations under consideration, in \Cref{sec.setup} we present the set-up of the benchmark problem under consideration, and in \Cref{sec.disctete-methods} we present details of the discretization methods under consideration. In \Cref{sec.numex}, we present the results achieved with the presented methods for the benchmark under consideration here and discuss these results. Finally, we give some concluding remarks in \Cref{sec:conclusions}.

\section{Setting the problem and model equations}
\label{sec.equations}

Let $(u,p)$ be a solution of the time-dependent Navier-Stokes
equations in a domain $\Omega\subset\mathbb{R}^d$ containing an obstacle with
boundary $\Gamma\subset\bo$. This fulfills the equations
\begin{subequations} \label{eqn:firstnavst}
\begin{align}
    && \partial_t u-\nu\Delta  u +  u\cdot\nabla u + \nabla p &= 0 &&\text{in}\;\Omega,\\
    && \sdiv u &=0&&\text{in}\;\Omega,
\end{align}
\end{subequations}
with the kinematic viscosity $\nu$, and together with boundary conditions
\begin{subequations}\label{eqn:bceesnavst}
\begin{align}
    u=g&\text{ on }\partial\Omega\backslash\Gamma\\
    u=0&\text{ on }\Gamma.
\end{align}
\end{subequations}
For the well-posedness of these equations, we refer to \cite{giraultraviart}.

\subsection{Weak formulation of the Navier-Stokes equations}
We 
consider a finite element discretizations of \eqref{eqn:firstnavst}.
To this end, we observe that the Navier-Stokes equations can be written in a weak
(or variational) form as follows: Find $(u,p) \in H^1_g(\Omega)\times L^2_0(\Omega)$,
such that
\begin{equation}\label{eqn:varform}
    (\partial_t u,v)_{L^2(\Omega)} + a(u, v) +c(u,u,v) + b(v, p) + b(u, q) = 0
\end{equation}
for all $(v,q)\in (H^1_0(\Omega))^d\times L^2_0(\Omega)$.
The space $H^1_g(\Omega)$ is the space of vector-valued $H^1$ functions with
trace $g$ on the boundary, and $H^1_0(\Omega)$ is the space of $H^1$
functions with trace $0$ on the boundary. The space $L^2_0(\Omega)$ is the space of
$L^2$ functions with mean zero. The bilinear forms $a(\cdot, \cdot)$ and
$b(\cdot, \cdot)$ are defined as
\begin{align*}
    a(u, v) \coloneqq \int_\Omega \frac{\nu}{2} \mathcal{D}(u):\mathcal{D}(w)\dif x
    \qquad\text{and}\qquad
    b(v, q) \coloneqq -\int_\Omega q (\nabla \cdot v) \dif x,
\end{align*}
respectively, where $\dv=\nabla v+\nabla v^t$ and the colon (:) indicates the
Frobenius inner-product of matrices. The convective (non-linear) term takes the form
\begin{equation}\label{eqn:conv}
    c(u,v,w)=c_\text{conv}(u, v, w) \coloneqq \int_\Omega (u\cdot\nabla v)\cdot w  \dif x .
\end{equation}
Other forms of the convective term may be found in the literature, see for example \cite{JohnBook16}.
For divergence-free $H^1$ functions, the various forms of the convective term
are equivalent. However, for non-divergence-free methods, other forms are beneficial and can be necessary for
non-$H^1$-conforming methods.

\subsection{Drag and lift evaluation}

The flow of fluid around an obstacle generates a force called drag, which is a
fundamental concept in fluid dynamics \cite{lamb1993hydrodynamics}.
It plays a critical role in determining the behavior of objects in flight and
has been studied since the time of d'Alembert \cite{lrsBIBjn}. The drag and
lift 
are given by the formula
\begin{equation}\label{draglift-direct}
    \beta(v) = \int_{\bog} \big((\nu\du - pI)n\big)\cdot v \dif s,
\end{equation}
where $I$ is the identity matrix, which consists of the viscous and pressure drag respectively. The drag is obtained using $v=(1,0)$ and the lift using $v=(0,1)$.
An alternative \cite{ref:refvalcyliftdragVolkerJohn} to evaluate $\beta$ is obtained through integration by parts, and
to test the weak formulation \eqref{eqn:varform} with a non-conforming test
function. This results in the functional
\begin{equation}\label{draglift-residual}
    \omega(v)= \int_\Omega (\partial_t u)\cdot v
     + \frac{\nu}{2}\du:\dv + (u\cdot\nabla u)\cdot v - p\sdiv v\dif x.
\end{equation}
Then, $\omega( v)=\beta( v)$ for all $ v\in H^1(\Omega)^d$
\cite{ref:refvalcyliftdragVolkerJohn,lrsBIBkj}.
This approach may be traced back to \cite{babuvska1984post}. It is also known
that, for the case of strongly imposed Dirichlet boundary conditions, this
approach of testing the residual with a non-conforming test-function doubles
the rate of convergence for the drag and lift, see for
example~\cite{braack2006solutions} for the proof in the steady case.

\section{Problem Set-up}
\label{sec.setup}
The considered set-up follows \cite{lrsBIBki,lrsBIBkl}.

\subsection{Geometry} 
We consider the domain
\begin{equation}\label{eqn:oneomega}
    \Omega=\{(x,y):-30<x<300,\,|y|<30,\;x^2+y^2>1\}.
\end{equation}
The boundary condition is set as $g = (1, 0)^T$ on the outer boundary and $g=(0,0)$
on the cylinder. The cylinder diameter is the reference length, so the Reynolds number
is given by $Re = 2 / \nu$. We consider the time interval $[0, 500]$.

\subsection{Quantities of interest}
For Reynolds numbers between 50 and 1000, the flow evolves into the von Karman vortex street \cite{Jac87,lrsBIBkl} with a periodic drag-lift profile. To measure the quality of the discrete solution, we consider whether a periodic solution is obtained and quantify this using the period of the drag-lift phase-diagram, i.e., the period for both drag and lift combined once the flow has fully developed.

\subsection{Initial condition}
For the initial condition, we consider the solution of the corresponding stationary Stokes problem. Using 0 as an initial condition and then ramping up the boundary condition would also make it possible to reach the desired state where the quantities of interest are evaluated.

\section{Discretization methods}
\label{sec.disctete-methods}

To discretize the weak formulation \eqref{eqn:varform}, we use finite element
methods. To this end, we take a simplicial mesh of the domain $\Omega$ with characteristic length
$h$, denoted as $\famlyT$. Finite element spaces then discretize
\eqref{eqn:varform} by approximating the spaces $H^1(\Omega)$ and
$L^2(\Omega)$ with spaces of piecewise polynomials on each element of the
mesh. To approximate the curved cylinder with higher order, the mesh on the cylinder is curved.

\subsection{Taylor-Hood}
As a base comparison, we consider an inf-sup stable (isoparametric) Taylor-Hood discretization of \eqref{eqn:varform}. 
It consists of $H^1$-conforming finite elements of order $k\geq 2$ for the
velocity and order $k-1$ for the pressure. We achieve optimal order of convergence of the finite element method by deforming the mesh with order $k$ to get the higher-order geometric approximation of the cylinder, see, e.g., \cite{brennerscott,EG21a}. This discretization will be denoted as $\THiso{k}$. When we use unaltered triangular elements to approximate the cylinder shape, we abbreviate the discretization as $\THk{k}$. Isoparametric Taylor-Hood elements are well-established and found to outperform non-deformed elements, even for the stationary flow around a cylinder \cite{JM01}. Nevertheless, to the best of our knowledge, the inf-sup stability of the isoparametric Taylor-Hood pair remains an open problem.

\subsubsection{Temporal Discretization}
To discretize in time, we consider both fully implicit and semi-implicit schemes. In the fully implicit case, we consider the backward-differentiation formula of order two or three for the time-derivative.

Let $\Delta t>0$ be a given time step, and $u_h^0$, $u_h^1$, $u_h^2$, be appropriate approximation of the initial condition and the solution at times $t=\Delta t$ and $t=2\Delta t$, respectively. In the fully implicit BDF3 scheme, we then solve the non-linear system
\begin{multline}\label{eqn:discrTH.bdf3}
    A_\text{bdf3}((u_h^n, p_h^n), (v_h,q_h)) \coloneqq  \int_{\Omega}\frac{11u_h^n - 18 u_h^{n-1} + 9 u_h^{n-2} - 2 u_h^{n-3}}{6\Delta t} v_h \dif x\\
     + a(u_h^n, v_h) + c(u_h^n,u_h^n,v_h) + b(v_h, p^n) + b(u_h^n, q_h) = 0
\end{multline}
in every time step $n\geq 3$. To get a scheme that is of order three in time, we initialize this scheme in the first two time step with the Crank-Nicolson ($\theta$-Scheme with $\theta=0.5$) time-discretization and 
\begin{multline}\label{eqn:discrTH.cn}
    A_\text{cn}((u_h^n, p_h^n), (v_h,q_h)) \coloneqq  \int_{\Omega}\frac{u_h^n - u_h^{n-1}}{\Delta t} v_h \dif x
     + \frac12a(u_h^n, v_h) + \frac12a(u_h^{n-1}, v_h) \\
     + \frac12c(u_h^n,u_h^n,v_h) + \frac12c(u_h^{n-1},u_h^{n-1},v_h) + b(v_h, p^n) + b(u_h^n, q_h) = 0.
\end{multline}
Note that we do not apply the $\theta$-Scheme to the pressure and divergence terms, to respect the fact that the pressure can be seen as the Lagrange multiplier for the divergence constraint. For an initial condition, we use the solution of the stationary Stokes problem with the same viscosity. Note that in our case, as we solve for $b(u_h^0, q_h)=0$, the above formulation is equivalent to solving for $p^{n-1/2}_n\approx p(t^n - 0.5 \Delta t)$ and enforcing $b((u_h^n + u_h^{n-1}) / 2, q_h) = 0$.

For the BDF2 scheme, the time derivative is approximated by $(\partial u_h, v_h)_{L^2(\Omega)}\approx \frac{1}{2\Delta t}(3 u_h^n - 4 u_h^{n-1} + u_h^{n-1}, v_h)_{L^2(\Omega)}$. We also initialize this scheme using the Crank-Nicolson scheme, although a single implicit Euler step would be sufficient to achieve a method of order two.

As a semi-implicit scheme, we consider the implicit-explicit (IMEX) based on the BDF2 time-discretization, also known as SBDF2 \cite{ARW95}. This uses the BDF2 time-discretization for the stiff Stokes part, while the non-linear convective term is treated with a second order extrapolation. The resulting scheme is
\begin{multline}\label{eqn:discrTH.sbdf2}
    \int_{\Omega}\frac{3u_h^n - 4 u_h^{n-1} + 2 u_h^{n-2}}{2\Delta t} v_h \dif x
     + a(u_h^n, v_h) + b(v_h, p^n) + b(u_h^n, q_h) \\
    + 2c(u_h^{n-1},u_h^{n-1},v_h) - c(u_h^{n-2},u_h^{n-2},v_h) = 0.
\end{multline}
Consequently, we need to solve the same \emph{linear} system in every time-step, so that we can reuse the matrix factorization in every time step. However, the explicit treatment of the convective term leads to a CFL condition. For $H^1$-conforming schemes for liner transport problems, this restriction can be as restrictive as $\Delta t \leq c(k) h^2$ \cite{EG21b}, depending on the explicit scheme used, while for discontinuous elements it can be shown that $\Delta t \leq chk^{-2}$ is needed for stability \cite{HW08}. Nevertheless, experimentally, it has been observed that $\Delta t \leq c h k^{-3/2}$ is sufficient in practice \cite{von2018implicit}.

\subsection{EMAC Formulation}
A common theme in the development of discretization methods is that schemes that include more physics lead to more accurate and stable solutions, especially over longer time intervals. While the above Taylor-Hood discretization is a prevalent method, the conservation properties of the continuous equation \eqref{eqn:varform}, such as those of energy, momentum and angular momentum, do not transfer to the discrete equation \eqref{eqn:discrTH.bdf3}. This is because $\nabla\cdot u_h \neq 0$ on the discrete level for Taylor-Hood elements. To recover these conservation properties without changing the finite element space, different formulations of the non-linear convective term are possible. For example, the skew-symmetric formulation $c_\text{skew}(u,u,v) = (u\cdot\nabla u + 0.5\diver(u)u, v)_{L^2(\Omega)}$ is energy conserving but does not conserve momentum and angular momentum, while $c_\text{cons}(u,u,v) = (\nabla\cdot(u\otimes u), v)_{L^2(\Omega)}$ conserves momentum and angular momentum but not energy. To conserve all three quantities, the so-called EMAC formulation \cite{CHOR17} of the non-linear term
\begin{equation}\label{eqn:conv-emac}
    c_\text{emac}( u, u, w) = (2(\du u), w) + ((\nabla \cdot  u) u, w),
\end{equation}
was developed. This form is consistent with the convective form $c_\text{conv}$ when $\nabla \cdot u=0$ but preserves energy, momentum and angular momentum even with $\nabla \cdot u \ne 0$. However, this scheme changes the pressure to $p_\text{kin} - \frac12 \vert u\vert^2$, where $p_\text{kin}$ is the kinematic pressure obtained with the standard convective term. If the forces on the cylinder are computed directly using \eqref{draglift-direct} rather than the residual formulation \eqref{draglift-residual}, we must take the change of pressure into account to obtain the correct result. The volumetric formulation has the additional advantage that no modification is needed.

We consider the same time-discretization schemes as for the standard Taylor-Hood method. This corresponds to \eqref{eqn:discrTH.bdf3}, \eqref{eqn:discrTH.cn} and \eqref{eqn:discrTH.sbdf2}, where the non-linear term $c(\cdot,\cdot,\cdot)$ is replaced with $c_\text{emac}(\cdot,\cdot,\cdot)$. We note that for the IMEX scheme, we can no longer expect the EMAC properties to hold due to the explicit treatment of the convective term. Furthermore, we note that IMEX methods using EMAC in conjunction with the exponential scalar auxiliary variable technique were recently analyzed in \cite{LLJ25}.

\subsection{Meshes}
We consider two different meshing approaches in two different software environment.

\subsubsection{Curved meshes}
We construct the mesh of the geometry \eqref{eqn:oneomega} using \texttt{Netgen}~\cite{schoeberl_netgen}. The domain is meshed over several levels such that we have a global mesh size $h_\text{max}= 8\cdot 2^{-\ell}$ for $\ell=0,1,2,3$, a smaller mesh size $h = h_\text{max}/2$ in a square of side-length 10, centered around the cylinder and a second local mesh size $h = h_{\text{max}} / 100$ on the surface of the cylinder. The resulting mesh is curved with the same order as the velocity finite element space, resulting in an isoparametric discretization. In what follows, we refer to this as \emph{mesh 1} and state $\ell$ to indicate the resolution.

\subsubsection{Straight meshes}
These meshes are constructed using \texttt{FreeFem++}\cite{freefem}.

We construct \emph{mesh 2}, splitting the left side into 20 segments, the first 80 units along the top and bottom into 35 segments, the remainder of the top and bottom into 25 segments, the right boundary into 12 segments, and the cylinder into 250  segments. The resulting mesh is refined globally once and then twice more in the box where $-15 < x < 70$ and $-15 < y < 15$. The $\THk{2}$ finite element space on this mesh has about $609\times 10^3$ degrees of freedom.

\emph{Mesh 3} is constructed similarly. Here, we divide the right sections along the top and bottom of the domain into 50 segments (rather than 25). We refine the resulting mesh globally and then once again in the box $-10 < x,y < 10$. The $\THk{3}$ finite element space on this mesh has about $632\times 10^3$ degrees of freedom.

\subsection{Non-linear Solver}

In every time step, we must solve the non-linear set of equations \eqref{eqn:discrTH.bdf3} or \eqref{eqn:discrTH.cn}. This may also be written as
\begin{equation*}
    \text{Find } U\in X\text{ such that } \quad A(U)(V) = F(V)\quad  \text{for all }V\in X,
\end{equation*}
where $A(\cdot)(\cdot)$ is linear in the second argument. With the Gâteaux derivative
\begin{equation*}
    A'(U)(W,V) \coloneqq \dod{}{s} A(U + sW)(V)\Big\vert_{s=0},
\end{equation*}
we consider a quasi-Newton method with a line-search globalization strategy in conjunction with the isoparametric Taylor-Hood elements implemented in \texttt{NGSolve}. In particular, if a full quasi-Newton step does not lead to a decrease in the residual, we consider a damped Newton update and decrease the damping parameter (i.e, a line-search) until the residual has decreased or a maximum number of search steps is reached. The procedure is summarized in the following.

\begin{algo}[Quasi-Newton with line-search globalization]~\\
\begin{algorithm}[H]
    \KwIn{Initial guess $U^0$, tolerance $\epsilon$, update tolerance $\delta<1$, damping factor $\gamma < 1$, maximum number line search steps $N$.}
    Compute factorization of Jacobian $D = A'(U^0)$\;
    Compute residual norm $\rho_0 = \Vert F - A(U^{0}) \Vert$\;
    i = 1\;
    \While{$\rho_{i-1} \geq \epsilon$}{
        Solve $D W^i = F - A(U^{i-1})$ for Newton update $W^i$\;
        $\omega_1 = 1$\;
        \For{$j=1,\dots, N$}{
            Update solution $U^{i,j} = U^{i-1} + \omega_j W^i$\;
            Compute residual norm $\rho_{i,j} = \Vert F - A(U^{i,j}) \Vert$\;
            \If{$\rho_{i,j} < \rho_{i-1}$}{
                Accept solution $U^i \to U^{i,j}$\;
                $\rho_i = \rho_{i,j}$\;
            }
            \Else{
                $\omega_{j+1} = \gamma \omega_j$\;
                $j\to j+1$\;
            }
        }
        \If{$j=N$}{
            Accept solution $U^i \to U^{i,N}$\;
            $\rho_i = \rho_{i,N}$\;
        }
        \If{$\rho_{i} > \delta \rho_{i-1}$}{
            Update factorization of Jacobian $D = A'(U^i)$.
        }
        $i \to i+1$\;
    }
    \KwResult{(Approximate) Solution $U^i$ of $A(U)=F$.}
\end{algorithm}
\end{algo}

In our implementation of Step 1, we reuse the last factorization of the Jacobian from the previous Newton solve, if available. For the initial guess, we take the solution from the previous time step. The residual tolerance is $\epsilon=10^{-11}$, the update tolerance is $\delta=0.1$, the damping factor is $\gamma=0.5$, and we consider a maximum number of $N=8$ line search steps. In particular, we found that if any damping is necessary in a quasi-Newton step, i.e, $j>1$, we also have that $\rho_{i} > \delta \rho_{i-1}$, so we update the Jacobian factorization for the next Newton step.

The (non-isoparametric) Taylor-Hood implementation using \texttt{FreeFem++} uses a standard full Newton scheme.

\section{Computational experiments}
\label{sec.numex}

We consider Reynolds numbers $500$ and $1000$ and compute the problem up to time $t=500$. We estimate the period of the drag-lift trajectory using a second-order approach described in \cite{lrsBIBkl}, between $t=280$ and $t=480$, during which the vortex street is fully developed. To compute the period, we need to truncate the time at the end of the computation, as our algorithm looks ahead, based on an approximation of the period, to find the closest point on the trajectory.

Running our computations, we found that the computations crashed in some cases. For the fully implicit methods, this occurred whenever the quasi-Newton method failed to converge, even after twenty Newton iteration steps, with up to 8 line search steps per iteration. Whenever the simulation failed before the vortex street was fully developed, we report the result as “solver failed” and the time at which the simulation failed. Sometimes, the solver failed after the vortex street was developed, allowing us to still approximate the vortex shedding frequency and resulting period. These cases are indicated by $t_\text{max} < 480$ in the tables below.

\subsection{Reynolds 500}

We present the results in \Cref{tab:TH500bdf3,tab:emac500bdf3,tab:emac500bdf2}. In \Cref{tab:TH500bdf3}, we see the results for the standard Taylor-Hood method, using the convective form \eqref{eqn:conv} of the non-linear term for $k=2,3,4$, with fully implicit BDF3 time-stepping. Here, we observe that the lowest-order case, $k=2$, yields poor results. While we get a periodic solution on the finest mesh with $k=2$, we have an error of $3.53\%$ in the Strouhal period. The higher-order results are significantly better. We get a periodic solution for coarser meshes, and we obtain the reference solution for $k=3$ on the finest mesh. Note that the latter is also the most resolved case, with approximately $427\times 10^3$ unknowns in the system after static condensation, cf. \Cref{tab:runtimes}.

\begin{table}
    \centering
    \caption{Results for $Re=500$ using the convective formulation and fully implicit BDF3 time stepping with $\Delta t=0.01$.}
    \label{tab:TH500bdf3}
    \begin{tabular}{llrccc}
        \toprule
        Discr. & Mesh & $t_\text{max}$ & Drag & Period & $\Vert\nabla\cdot u\Vert_{\ell^\infty(L^2(\Omega))}$\\
        \midrule
        $\THiso{2}$ & 1, $\ell=0$ & 480.00 & 0.870 & $\phantom{1}9.736\pm2.62830$  & $9.5 \cdot 10^{0}$ \\
        $\THiso{2}$ & 1, $\ell=1$ & 480.00 & 1.140 & $10.997\pm2.00802$ & $1.8 \cdot 10^{1}$ \\
        $\THiso{2}$ & 1, $\ell=2$ & 480.00 & 1.178 & $\phantom{1}9.834\pm1.31278$  & $8.2 \cdot 10^{0}$ \\
        $\THiso{2}$ & 1, $\ell=3$ & 480.00 & 1.400 & $\phantom{1}9.131\pm0.00005$  & $4.4 \cdot 10^{0}$ \\[3pt]
        $\THiso{3}$ & 1, $\ell=0$ & 480.00 & 1.471 & $\phantom{1}8.887\pm0.04725$  & $1.1 \cdot 10^{1}$ \\
        $\THiso{3}$ & 1, $\ell=1$ & 480.00 & 1.432 & $\phantom{1}8.934\pm0.00269$  & $1.2 \cdot 10^{1}$ \\
        $\THiso{3}$ & 1, $\ell=2$ & 480.00 & 1.434 & $\phantom{1}8.752\pm0.00007$  & $9.2 \cdot 10^{0}$ \\
        $\THiso{3}$ & 1, $\ell=3$ & 480.00 & 1.445 & $\phantom{1}8.823\pm0.00211$  & $1.0 \cdot 10^{1}$ \\[3pt]
        $\THiso{4}$ & 1, $\ell=0$ & 480.00 & 1.438 & $\phantom{1}8.834\pm0.01065$  & $1.2 \cdot 10^{1}$ \\
        $\THiso{4}$ & 1, $\ell=1$ & 480.00 & 1.441 & $\phantom{1}8.883\pm0.00007$  & $1.3 \cdot 10^{1}$ \\
        $\THiso{4}$ & 1, $\ell=2$ & 480.00 & 1.441 & $\phantom{1}8.836\pm0.00430$  & $9.7 \cdot 10^{0}$ \\
        \cmidrule(lr){1-6}
        \multicolumn{3}{r}{Reference~\cite{lrsBIBkl}:} & 1.448 & \multicolumn{1}{l}{\phantom{1}8.82} \\
        \bottomrule
   \end{tabular}
\end{table}

In \Cref{tab:emac500bdf3}, we see the results obtained using the EMAC convective term \eqref{eqn:conv-emac} for the same meshes, polynomial orders, and time-stepping scheme as those used above. For the lowest-order case, $k=2$, we observe that a periodic solution is obtained on coarser meshes. Furthermore, the resulting Strouhal period on the finest mesh is about twice as accurate as before, with a relative error of $1.79\%$. However, problems appear in the case of higher-order isoparametric elements, and none of the simulations reached the final time $t=500$. Nevertheless, for the cases where the vortex street reached the fully developed state, the results are again accurate. We note that these cases are on fine meshes using higher-order elements.

\begin{figure}
    \centering
    \includegraphics{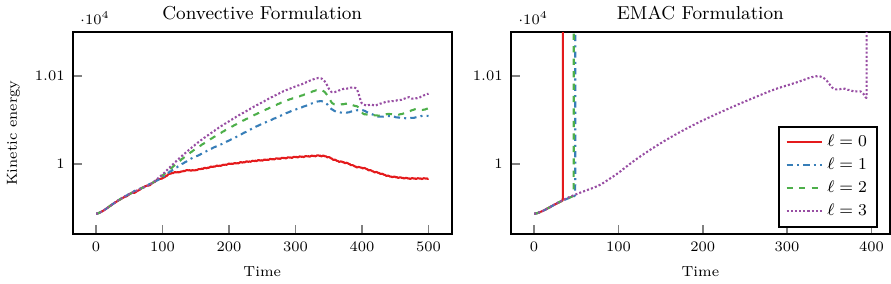}
    \caption{Kinetic energy over time for $k=3$ using fully implicit BDF3 time-stepping.}
    \label{fig:Re500-energy}
\end{figure}

\begin{figure}
    \centering
    \includegraphics{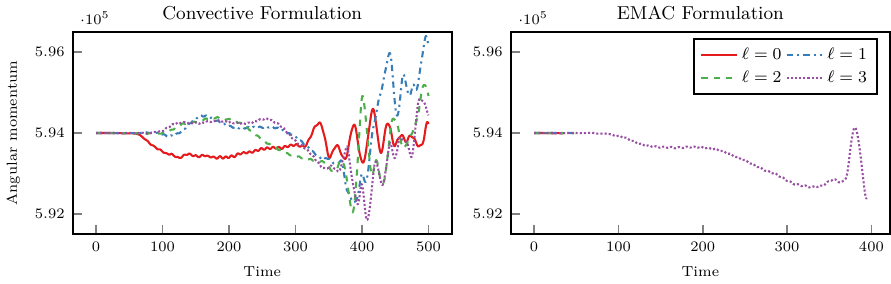}
    \caption{Angular momentum over time for $k=3$ using fully implicit BDF3 time-stepping.}
    \label{fig:Re500-angmom}
\end{figure}

\begin{figure}
    \centering
    \includegraphics{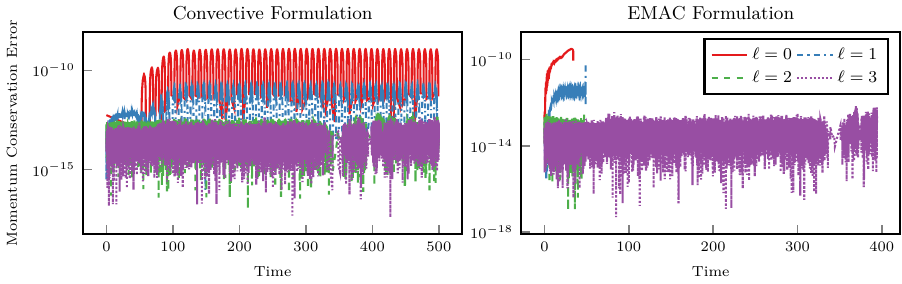}
    \caption{Momentum conservation error over time for $k=3$ using fully implicit BDF3 time-stepping.}
    \label{fig:Re500-mom}
\end{figure}

Examining the results from the finer meshes with polygonal approximation of the cylinder, we observe that the reference values are achieved. In particular, the case $k=3$ did not crash. However, we note that we did observe crashes in the \texttt{FreeFem++} implementation when similar coarse meshes to those used with \texttt{NGSolve} were used. Consequently, it appears that sufficient mesh resolution is more critical for a stable computation using the EMAC formulation with higher-order elements. Looking at the development of kinetic energy, angular and linear momentum for the case $k=3$ in \Cref{fig:Re500-energy,fig:Re500-angmom,fig:Re500-mom}, respectively, we see that the crashes (failure for Newton's method to solve the resulting linear system) are characterized by an explosion in the kinetic energy of the solution.

\begin{table}
    \centering
    \caption{Results for $Re=500$ using the EMAC formulation and fully implicit BDF3 time stepping.}
    \label{tab:emac500bdf3}
    \begin{tabular}{lllrccc}
        \toprule
        Discr. & Mesh & $\Delta t$ & $t_\text{max}$ & Drag & Period & $\Vert\nabla\cdot u\Vert_{\ell^\infty(L^2(\Omega))}$\\
        \midrule
        $\THiso{2}$ & 1, $\ell=0$ & 0.01 & 480.00 & 1.073\phantom{0}  & $11.381\pm2.05767$ & $1.0\cdot 10^{1}$ \\
        $\THiso{2}$ & 1, $\ell=1$ & 0.01 & 480.00 & 1.321\phantom{0}  & $\phantom{1}9.897\pm0.12250$ & $1.0\cdot 10^{1}$ \\
        $\THiso{2}$ & 1, $\ell=2$ & 0.01 & 480.00 & 1.372\phantom{0}  & $\phantom{1}9.838\pm0.00316$ & $8.7\cdot 10^{0}$ \\
        $\THiso{2}$ & 1, $\ell=3$ & 0.01 & 480.00 & 1.381\phantom{0}  & $\phantom{1}8.978\pm0.00023$ & $6.2\cdot 10^{0}$ \\[3pt]
        $\THk{2}$   & 2           & 0.02 & 480.00 & 1.4470 & $\phantom{1}8.819\pm 0.0044\phantom{0}$\\[3pt]
        $\THiso{3}$ & 1, $\ell=0$ & 0.01 & 34.17 && Solver failed\\
        $\THiso{3}$ & 1, $\ell=1$ & 0.01 & 48.86 && Solver failed\\
        $\THiso{3}$ & 1, $\ell=2$ & 0.01 & 46.99 && Solver failed\\
        $\THiso{3}$ & 1, $\ell=3$ & 0.01 & 374.30 & 1.445\phantom{0}  & $\phantom{1}8.820\pm0.00080$ & $1.0\cdot 10^{1}$ \\[3pt]
        $\THk{3}$   & 3           & 0.02 & 480.00 & 1.4472 & $\phantom{1}8.820\pm 0.0050\phantom{0}$ \\[3pt]
        $\THiso{4}$ & 1, $\ell=0$ & 0.01 & 57.48 && Solver failed\\
        $\THiso{4}$ & 1, $\ell=1$ & 0.01 & 397.20 & 1.457\phantom{0} & $\phantom{1}8.799\pm0.00156$ & $1.3\cdot 10^{1}$ \\
        $\THiso{4}$ & 1, $\ell=2$ & 0.01 & 380.93 & 1.435\phantom{0} & $\phantom{1}8.881\pm0.00048$ & $1.2\cdot 10^{1}$ \\
        \cmidrule(lr){1-7}
        \multicolumn{4}{r}{Reference~\cite{lrsBIBkl}:} & 1.448\phantom{0} & \multicolumn{1}{l}{\phantom{1}8.82} \\
        \bottomrule
   \end{tabular}
\end{table}

To check whether the above behavior of the EMAC scheme can be explained by the near instability of the BDF3 scheme, the non-linear solver, or the choice of time-step, we compute the case $k=3$ using the more stable BDF2 scheme in time, the semi-implicit SBDF2 scheme and with a smaller time-step. The results are summarized in \Cref{tab:emac500bdf2}. Here, we see that computations fail at a very similar point in time, independent of the time-stepping scheme and time-step size for a given spatial discretization. We conclude that the previously observed failures are indeed due to the spatial rather than the temporal discretization.

\begin{table}
    \centering
    \caption{Results for $Re=500$ using the EMAC formulation with different time steps and time-stepping schemes.}
    \label{tab:emac500bdf2}
    \begin{tabular}{llllrccc}
        \toprule
        Discr. & Mesh & Scheme & $\Delta t$ & $t_\text{max}$ & Drag & Period & $\Vert\nabla\cdot u\Vert_{\ell^\infty(L^2(\Omega))}$\\
        \midrule
        $\THiso{3}$ & 1, $\ell=0$ & BDF2 & 0.01 & 34.16 && Solver failed\\
        $\THiso{3}$ & 1, $\ell=1$ & BDF2 & 0.01 & 48.85 && Solver failed\\
        $\THiso{3}$ & 1, $\ell=2$ & BDF2 & 0.01 & 46.99 && Solver failed\\
        $\THiso{3}$ & 1, $\ell=3$ & BDF2 & 0.01 & 374.45 & 1.444 & $8.820\pm0.00079$ & $1.0\cdot 10^{1}$ \\[3pt]
        $\THiso{3}$ & 1, $\ell=0$ & SBDF2 & 0.003 & 34.25 && Solver failed\\
        $\THiso{3}$ & 1, $\ell=1$ & SBDF2 & 0.0015 & 48.91 && Solver failed\\
        $\THiso{3}$ & 1, $\ell=2$ & SBDF2 & 0.00075 & 47.08 && Solver failed\\[3pt]
        $\THiso{3}$ & 1, $\ell=0$ & BDF3 & 0.005 & 34.18 && Solver failed\\
        $\THiso{3}$ & 1, $\ell=1$ & BDF3 & 0.005 & 48.86 && Solver failed\\
        $\THiso{3}$ & 1, $\ell=2$ & BDF3 & 0.005 & 47.00 && Solver failed\\
        $\THiso{3}$ & 1, $\ell=3$ & BDF3 & 0.005 & 374.30 & 1.445 & $8.820\pm0.00108$ & $1.0\cdot 10^{1}$ \\
        \cmidrule(lr){1-8}
        \multicolumn{5}{r}{Reference~\cite{lrsBIBkl}:} & 1.448 & \multicolumn{1}{l}{8.82} \\
        \bottomrule
   \end{tabular}
\end{table}

We compare some basic computational quantities of our simulations in \Cref{tab:runtimes} for the different spatial and temporal discretizations. The computations were performed using four shared-memory parallel cores of an Intel Xeon Platinum 8360Y processor at 2.4 GHz on the University of Jena's scientific cluster, ``draco''. The resulting linear systems were solved using the direct solver \texttt{pardiso} via the Intel MKL. As we are using higher-order $H^1$-conforming elements, element internal degrees of freedom can be removed from the resulting linear systems using static condensation \cite{EG21b}. This can be done both for the IMEX scheme and for the Jacobian in the quasi-Newton method. Consequently, we report the full number of degrees of freedom (\texttt{ndof}) of the finite element space and the number of unconstrained remaining after static condensation (\texttt{cond}).

Looking at the results in \Cref{tab:runtimes}, we first note that the coarsest mesh leads to a small problem with only about $17\times10^3$ unknowns, while $k=3$ on mesh 1 with $\ell=3$ leads to the largest number of unknowns, both before and after static condensation. Examining the total number of quasi-Newton steps, we observe that the EMAC formulation typically necessitates more Newton iterations. Furthermore, we also observe that the EMAC formulation requires significantly more updates of the Jacobian factorization, illustrating the more difficult non-linearity introduced by EMAC. More surprisingly, we see that some runs using the convective formulation required only a few updates of the Jacobian during the approximately $2.5\times10^5$ Newton steps, averaging to 5 iterations per time step. This makes the approach highly competitive against IMEX schemes, as the fully implicit scheme does not require fulfillment of a CFL condition for stability.

\begin{table}
    \centering
    \caption{Comparison of computational run times for $Re=500$ for different schemes. The total number of degrees of freedom are abbreviated by \texttt{ndof}, and the number of unconstrained dofs after static condensation by \texttt{cond}.}
    \label{tab:runtimes}
    \begin{tabular}{llllrrrrr}
        \toprule
        Scheme & Discr. & Mesh & \texttt{ndof}(\texttt{cond}) $\times 10^3$ & $\Delta t$ & Newton Steps & Jac. Upd. &  Time (s)\\
        \midrule
        BDF3 conv. & $\THiso{2}$ & 1, $\ell=1$ & 18(17)   & 0.01 & 248150 &  44 &   6019 \\ 
        BDF3 EMAC  & $\THiso{2}$ & 1, $\ell=1$ & 18(17)   & 0.01 & 271135 & 507 &   5774 \\[2pt] 
        BDF3 conv. & $\THiso{2}$ & 1, $\ell=2$ & 60(58)   & 0.01 & 254779 &   4 &  17006 \\ 
        BDF3 EMAC  & $\THiso{2}$ & 1, $\ell=2$ & 60(58)   & 0.01 & 263836 & 598 &  17657 \\[2pt] 
        BDF3 conv. & $\THiso{2}$ & 1, $\ell=3$ & 218(213) & 0.01 & 271480 & 186 &  78858 \\ 
        BDF3 EMAC  & $\THiso{2}$ & 1, $\ell=3$ & 218(213) & 0.01 & 279464 & 711 &  86270 \\[2pt] 
        BDF3 conv. & $\THiso{3}$ & 1, $\ell=2$ & 146(116) & 0.01 & 245878 & 317 &  50549 \\ 
        BDF3 conv. & $\THiso{3}$ & 1, $\ell=3$ & 530(427) & 0.01 & 282482 &   3 & 212431 \\ 
        BDF3 conv. & $\THiso{4}$ & 1, $\ell=2$ & 270(174) & 0.01 & 244402 & 312 &  80875 \\ 
        \bottomrule
    \end{tabular}
\end{table}

\subsection{Reynolds 1000}
The results are summarized in \Cref{tab:TH1000bdf3,tab:emac1000bdf3}. In \Cref{tab:TH1000bdf3}, we now see the convective formulation also begins to struggle across all considered polynomial orders. However, we observe that the crash is delayed by both increasing the mesh resolution and the polynomial order. This contrasts with the previous EMAC computations, where higher-order terms did not necessarily improve the robustness. On the finest discretization ($\ell=3, k=3$), we see that the reference quantities are almost obtained.

Looking at the results for the EMAC formulation in Table 6, we make several key observations. First, we observe that for the lowest-order case, $k = 2$, the method appears more stable than the convective formulation, with only the coarsest discretization failing. Using higher-order isoparametric elements, we observe that the simulations fail earlier than in the case of $Re = 500$, highlighting the increased difficulty associated with higher Reynolds numbers. Nevertheless, the finest discretizations (using straight elements) manage to resolve the flow, and we achieve close to the reference results using P3/P2 elements on mesh 3, illustrating again that sufficient resolution is essential here.

\begin{table}
    \centering
    \caption{Results for $Re=1000$ using the convective formulation and fully implicit BDF3 time stepping with $\Delta t=0.01$.}
    \label{tab:TH1000bdf3}
    \begin{tabular}{llrccc}
        \toprule
        Discr. & Mesh & $t_\text{max}$ & Drag & Period & $\Vert\nabla\cdot u\Vert_{\ell^\infty(L^2(\Omega))}$\\
        \midrule
        $\THiso{2}$ & 1, $\ell=0$ & 74.99 && Solver failed\\
        $\THiso{2}$ & 1, $\ell=1$ & 63.31 && Solver failed\\
        $\THiso{2}$ & 1, $\ell=2$ & 95.06 && Solver failed\\
        $\THiso{2}$ & 1, $\ell=3$ & 480.00 & 1.402 & $8.621\pm0.01750$ & $1.2\cdot 10^{1}$ \\[3pt]
        $\THiso{3}$ & 1, $\ell=0$ & 94.40 && Solver failed\\
        $\THiso{3}$ & 1, $\ell=1$ & 51.63 && Solver failed\\
        $\THiso{3}$ & 1, $\ell=2$ & 339.35 & 1.533 & $8.367\pm0.00119$ & $9.7\cdot 10^{0}$ \\
        $\THiso{3}$ & 1, $\ell=3$ & 379.49 & 1.534 & $8.357\pm0.00063$ & $1.1\cdot 10^{1}$ \\[3pt]
        $\THiso{4}$ & 1, $\ell=0$ & 54.74 && Solver failed\\
        $\THiso{4}$ & 1, $\ell=1$ & 480.00 & 1.526 & $8.432\pm0.45562$ & $1.8\cdot 10^{1}$ \\
        $\THiso{4}$ & 1, $\ell=2$ & 480.00 & 1.529 & $8.371\pm0.03071$ & $1.9\cdot 10^{1}$ \\
        \cmidrule(lr){1-6}
        \multicolumn{3}{r}{Reference~\cite{lrsBIBkl}:} & 1.54\phantom{9} & \multicolumn{1}{l}{8.36} \\
        \bottomrule
   \end{tabular}
\end{table}

\begin{table}
    \centering
    \caption{Results for $Re=1000$ using EMAC Taylor-Hood and fully implicit BDF3 time stepping with $\Delta t=0.01$.}
    \label{tab:emac1000bdf3}
    \begin{tabular}{lllrccc}
        \toprule
        Discr. & Mesh & $\Delta t$ & $t_\text{max}$ & Drag & Period & $\Vert\nabla\cdot u\Vert_{\ell^\infty(L^2(\Omega))}$\\
        \midrule
        $\THiso{2}$ & 1, $\ell=0$ & 0.01 & 99.11 && Solver failed\\
        $\THiso{2}$ & 1, $\ell=1$ & 0.01 & 480.00 & 1.332\phantom{0} & $9.443\pm1.23356$ & $1.6\cdot 10^{1}$ \\
        $\THiso{2}$ & 1, $\ell=2$ & 0.01 & 480.00 & 1.405\phantom{0} & $9.408\pm0.05105$ & $1.6\cdot 10^{1}$ \\
        $\THiso{2}$ & 1, $\ell=3$ & 0.01 & 480.00 & 1.421\phantom{0} & $8.775\pm0.00059$ & $1.4\cdot 10^{1}$ \\[3pt]
        $\THk{2}$   & 2           & 0.02 & 480.00 & 1.5678 & $8.280\pm 0.0447\phantom{0}$ &\\[3pt]
        $\THiso{3}$ & 1, $\ell=0$ & 0.01 & 20.00 && Solver failed\\
        $\THiso{3}$ & 1, $\ell=1$ & 0.01 & 22.42 && Solver failed\\
        $\THiso{3}$ & 1, $\ell=2$ & 0.01 & 14.90 && Solver failed\\
        $\THiso{3}$ & 1, $\ell=3$ & 0.01 & 31.28 && Solver failed\\[3pt]
        $\THk{3}$   & 3           & 0.02 & 480.00 & 1.5413 & $8.371\pm 0.0101\phantom{0}$ \\[3pt]
        $\THiso{4}$ & 1, $\ell=0$ & 0.01 & 38.94 && Solver failed\\
        $\THiso{4}$ & 1, $\ell=1$ & 0.01 & 19.93 && Solver failed\\
        $\THiso{4}$ & 1, $\ell=2$ & 0.01 & 38.73 && Solver failed\\
        \cmidrule(lr){1-7}
        \multicolumn{4}{r}{Reference~\cite{lrsBIBkl}:} & 1.54\phantom{00} & \multicolumn{1}{l}{8.36}\\
        \bottomrule
   \end{tabular}
\end{table}

\section{Conclusions}
\label{sec:conclusions}

We presented a computational study of a stress test simulation of the flow around a cylinder in a large domain over a long period at Reynolds numbers 500 and 1000. This is a challenging set-up due to the long-time accuracy required to obtain an accurate solution to this problem. The presented methods are based on prevalent Taylor-Hood elements, using both straight and isoparametric elements. In conjunction with these spaces, we considered the Navier-Stokes equations, both in the standard convective form of the non-linear term and in the EMAC formulation. The latter is designed such that energy, momentum, and angular momentum are conserved, even when $\diver(u_h)\neq0$, which is the case for the Taylor-Hood finite element pair. Including these physical properties in the discrete scheme promises a more accurate scheme, specifically for the long simulation required here. For the temporal discretization, we considered the fully implicit BDF3 and BDF2 schemes, as well as the implicit-explicit SBDF2 scheme.

For the lowest-order Taylor-Hood pair, P2/P1, we found that the EMAC scheme indeed yields better results for the same spatial discretization. However, the increased non-linearity required more frequent of the Jacobian in our quasi-Newton scheme, such that the time to solution was longer. Nevertheless, for $k = 2$ and $\ell = 3$, the Strouhal period error almost halved from 3.53\% to 1.79\% by changing from the convective form to the EMAC formulation while needing only approximately 10\% more time-to-solution.

However, for higher-order schemes, $k\geq3$, we found that the EMAC formulation is more sensitive to coarse meshes, and the simulations crashed due to a failure of Newton's scheme that coincided with an exponential energy growth. Nevertheless, on finer meshes, the EMAC formulation obtained the reference results from the literature, even at $Re = 1000$. Furthermore, the EMAC formulation yielded a stable simulation at $Re = 1000$ for $k = 2$, whereas the convective formulation crashed on the first three mesh levels.

\section*{Data Availability Statement}
Both the \texttt{NGSolve} and \texttt{FreeFem++} scripts used to realize the presented results is available on github and archived on zenodo at \url{https://doi.org/10.5281/zenodo.15869558}.

\section*{Acknowledgments}
Author LGR acknowledges partial support from the US National Science Foundation grant DMS 2152623.

\printbibliography

\end{document}